\documentclass{amsart}

\usepackage{amsmath}
\usepackage{amscd}
\usepackage{subfigure}

\usepackage{epsfig}

\newtheorem{theorem}{Theorem}[section]
\newtheorem{corollary}[theorem]{Corollary}

\newtheorem{problem}[theorem]{Problem}

{
\theoremstyle{plain}
\newtheorem{definition}[theorem]{Definition}

\newtheorem{ack}{Acknowledgements}
}

\newcommand{\A}{\ensuremath{\mathcal A}}
\renewcommand{\S}{\ensuremath{\mathcal S}}
\renewcommand{\P}{\ensuremath{\mathbb P}}
\newcommand{\R}{\ensuremath{\mathcal R}}
\newcommand{\B}{\ensuremath{\mathcal B}}
\newcommand{\G}{\ensuremath{\mathcal G}}
\newcommand{\C}{\ensuremath{\mathbb{C}}}
\newcommand{\Q}{\ensuremath{\mathbb{Q}}}
\newcommand{\Z}{\ensuremath{\mathbb{Z}}}
\newcommand{\rk}{\operatorname{rk}}
\renewcommand{\Re}{\ensuremath{\mathbb{R}}}
\newcommand{\Aut}{\operatorname{Aut}}
\newcommand{\Ext}{\operatorname{Ext}}

\begin{document}
\title[Homotopy theory, II]{On the homotopy theory of arrangements, 
II}
\author[M. Falk]{Michael Falk}
\address{Northern Arizona University\\
Flagstaff, Arizona USA}
\email{Michael.Falk@nau.edu}
\author[R. Randell]{Richard Randell}
\address{University of Iowa\\
Iowa City, Iowa 52242 USA}
\email{randell@math.uiowa.edu}
\date{November 30, 1998}
\thanks{This paper is in final form and no version of it will be 
submitted for
publication elsewhere.}
\begin{abstract}
In ``On the homotopy theory of arrangements'' published in 1986 the 
authors
gave a comprehensive survey of the subject. This article updates and
continues the earlier article, noting some key open problems.
\end{abstract}

\maketitle

Let $M$ be the complement of a complex arrangement. Our interest here 
is in
the topology, and especially the homotopy theory of $M$, which turns 
out to
have a rich structure. In the first paper of this name \cite{FR1}, we assembled 
many of
the known results; in this paper we wish to summarize progress in the
intervening years, to reiterate a few key unsolved questions, and 
propose some new problems we find of interest.

In the first section we establish some terminology and notation, and 
discuss general homotopy classification problems. We introduce the
matroid-theoretic terminology that has become more prevalent in the 
subject in recent years. In this section we also sketch Rybnikov's 
construction of arrangements with the same matroidal structure but 
non-isomorphic fundamental groups.
In Section 2 we consider some algebraic properties 
of
the fundamental group of the arrangement. Properties of interest 
include the
lower central series, the Chen groups, the rational homotopy theory 
of the
complement, and the cohomology of the group. At the time of our first 
paper
many questions in this area were in flux, so we make a special effort 
here
to clarify the situation. The group cohomology is naturally of 
interest in
the third section as well, which focuses on when or if the complement 
is
aspherical. It is this property which fostered much of the initial
interest in arrangements (in the guise of the pure braid space); it 
is of
interest that the determination of when the complement is aspherical 
is far
from settled. Finally, in the fourth section we
consider what one might call the topology of the fundamental group.
We describe group presentations that have been discovered since the 
publication of \cite{FR1}, including the recent development of braided 
wiring diagrams. We also sketch the considerable progress in the 
study of the Milnor fiber associated with an arrangement.

In 1992 the long-awaited book {\em Arrangements of Hyperplanes,} by 
Peter
Orlik and Hiroaki Terao appeared, to the delight of all of us working 
in
arrangements. We refer the reader to this text as a general reference 
on arrangements, and adopt their notation and terminology 
except where specified. We also mention that perhaps the most 
interesting 
development in
arrangements in the last ten years involves the deep and fascinating
connections with hypergeometric functions. We are pleased to refer the
reader to the lecture notes of Orlik and Terao \cite{O3} from the 1998 Tokyo 
meeting for a 
comprehensive
exposition of this material.

\begin{section}{Combinatorial and topological structure}

One significant change in the study of the homotopy theory of arrangements since 
the publication of \cite{FR1} has been the introduction of 
matroid-theoretic terminology and techniques into the subject. In 
this section we review this approach and describe progress toward the 
topological classification of hyperplane complements. Refer to 
\cite{Wh1,Ox} 
for further details on matroids.

\begin{subsection}{The matroid of an arrangement}
Let $V=\C^\ell$ and let $\A=\{H_1, \ldots, H_n\}$ be a central 
arrangement of hyperplanes in $V$. For each hyperplane $H_i$ choose a 
linear form $\alpha_i\in V^*$ with $H_i=\ker(\alpha_i)$. The product 
$Q(\A)=\prod_{i=1}^n \alpha_i$ is the {\em defining polynomial} of 
the arrangement.

The {\em underlying matroid} $G(\A)$  of \A\ is by definition the 
collection of
subsets of $[n]:=\{1,\ldots, n\}$ given by $$G(\A)=\{S\subseteq [n] \ 
| \ \{\alpha_i \ | \ i \in S\} \ \text{is linearly dependent}\}.$$

Elements of $G=G(\A)$ are called {\em dependent sets}. Minimal 
dependent sets are called {\em circuits}. Independent sets and bases 
are defined in the obvious way. The {\em rank} $\rk(S)$ of a set $S 
\subseteq [n]$ is the size of a maximal independent subset of $S$. 
The rank of $G$ (or \A) is $\rk([n])$.
The {\em closure} $\overline{S}$ of a set $S$ is defined by 
$$\overline{S}=\bigcup\{T \subseteq [n] \ | \ T \supseteq S \ 
\text{and} \ \rk(T)=\rk(S)\}.$$

A set $S$ is {\em closed} if $\overline{S}=S$. Closed sets are also 
called {\em flats}. The collection of closed sets, ordered by 
inclusion, forms a geometric lattice $L(G)$ which is isomorphic to 
the intersection lattice $L(\A)$ defined and studied in \cite{OT}. 
The isomorphism $L(G)\to L(\A)$ is given by $S \mapsto \bigcap_{i \in 
S} H_i$.

Thus the matroid $G(\A)$ contains the same information as the 
intersection lattice $L(\A)$. One of the simple advantages of the 
matroid-theoretic approach is the fact that the matroid $G(\A)$ is 
determined uniquely by any of a number of different pieces of data 
besides the set of flats. For instance, the set of circuits, the rank 
function, or the set of bases, each determine the matroid, and thus 
the intersection lattice. Besides giving a nice conceptual framework 
for the combinatorial structure of arrangements, techniques and deep 
results from the matroid theory literature have been applied with 
some benefit in the study of  the topology of arrangements. 

The line generated by $\alpha_i$ in $V^*$ depends only on $H_i$, and 
thus \A\ determines a unique point configuration $\A^*$ in the 
projective 
space $\P(V^*)\cong \C P^{\ell-1}$. 
The dual point configuration $\A^*$ can be used to depict the 
combinatorial structure of an arrangement in case $\rk(\A)\leq 4$ if 
the defining forms $\alpha_i$ have real coefficients. (In this case 
\A\ is called a {\em complexified} arrangement.) One merely plots the 
points $\alpha_i$ in a
suitably chosen affine chart $\Re^{\ell-1}$ in the real projective 
space $\Re P^{\ell-1}$, for instance by scaling the $\alpha_i$ so 
that the coefficient of $x_1$ in each is equal to 1, and then 
ignoring this coefficient. Dependent flats of rank two (or three) are 
seen in these affine configurations as lines (or planes) containing 
more 
than two (or three) points. These lines and planes are usually 
explicitly indicated in the picture. This is especially useful for 
arrangements of rank four. Since the hyperplanes are indicated by 
points in $\Re^3$, they don't obscure the internal structure as a 
collection of affine planes in $\Re^3$ would (see Figure \ref{yuzfig}).
These depictions of projective point configurations are generalized 
to give {\em affine diagrams} of arbitrary matroids. Dependent flats 
are again explicitly indicated with ``lines" or ``planes," which in 
the general case may not be straight or flat in the euclidean 
sense. It is common to refer to flats of rank one, two, or three in an arbitrary 
matroid as  points, lines, or planes respectively. These diagrams are useful for 
the study of arrangements which are not complexified real arrangements (see 
Figures \ref{mac} and \ref{ryb}).

\end{subsection} 

\begin{subsection}{Basic topological results}
\label{top}
The seminal result in the homotopy theory of arrangements is the 
calculation of the cohomology algebra of the complement 
$M=M(\A):=\C^\ell - \bigcup_{i=1}^n H_i$ by 
Orlik and Solomon \cite{OS1}. Motivated by work of Arnol'd \cite{A}, 
and using tools established by Brieskorn \cite{Br}, they gave a 
presentation 
of $H^*(M)$ in terms of generators and relations. The presentation 
$A(\A)$ depends only on the underlying matroid $G=G(\A)$, and 
is now called the 
{\em Orlik-Solomon (or $OS$) algebra} of $G$. Henceforth we will 
refer to the $OS$ algebra $A(\A)$ rather than the cohomology ring 
$H^*(M)$. The algebra $A(\A)$ is defined as the quotient of the 
exterior algebra on generators $e_1,\ldots, e_n$ by the ideal $I$ 
generated by ``boundaries" of dependent sets of $G$. See \cite{OT} 
for a precise definition.

This result of \cite{OS1} gave rise to a collection of ``homotopy type" 
conjectures, which assert that various homotopy invariants of the 
complement depend only on $G(\A)$.
A great deal of research in the homotopy theory of arrangements has 
been focused on conjectures of this type. Note that such conjectures 
may have ``weak'' or ``strong'' solutions: one may show that the 
invariant depends only on the matroid, or one may give an algorithm to 
compute the invariant from matroidal data.

The major positive result 
in this direction is the lattice-isotopy theorem, proved by the 
second author in \cite{R1}. It asserts that the homotopy type, indeed 
the diffeomorphism type of the complement remains constant through a 
``lattice-isotopy," that is, a one-parameter family of arrangements 
in which the intersection lattice, or equivalently, the underlying 
matroid remains constant. 

This result is often recast in terms of matroid realization 
spaces, which are related to the well-known ``matroid stratification" of the 
Grassmannian. We describe this connection. The defining forms 
$\alpha_i$ of \A\ can be identified with row vectors, and thus the 
arrangement \A\ can be identified with an $n \times \ell$ matrix $R$ 
over \C. This matrix is called a {\em realization} of the underlying 
matroid. Two realizations $R$ and $R'$ are equivalent if there is a 
nonsingular diagonal $n\times n$ matrix $S$ and a nonsingular 
$\ell\times\ell$ 
matrix $T$ such that $R'=SRT$. The corresponding arrangements will 
then be linearly isomorphic. The set of equivalence classes of 
realizations of a fixed matroid $G$ is called the (projective) 
realization space 
$\R(G)$ of $G$. Now assume the matrix $R$ has rank $\ell$, i.e., that 
\A\ is an essential arrangement. Then the column space of $R$ is an 
$\ell$-plane 
$P_R$ (sometimes denoted $P_\A$) in $\C^n$.
Note that an isomorphic copy of the arrangement 
$\A$ inside $P_R$ is formed by the intersection of $P_R$ with the 
coordinate hyperplanes in $\C^n$. Postmultiplying $A$ by a 
nonsingular 
matrix doesn't affect $P_R$. Thus we see that the realization space 
$\R(G)$ can be identified with a subset $\Gamma(G)$ of the space of 
orbits of the diagonal $(\C^*)^n$ action on the Grassmanian 
$\G_\ell(\C^n)$ of $\ell$-planes in $\C^n$. The subsets 
$\widehat{\Gamma}(G)=\{ P_R \ | \ R \ \text{is a realization of} \ 
G\}\subseteq \G_\ell(\C^n)$ 
are called
{\em matroid strata}, although they do not comprise a stratification 
in the usual sense, since the closure of a stratum 
may not be a union of strata \cite{Sturm1}. These strata play a central role in 
the 
theory of generalized hypergeometric functions, especially when the 
original arrangement \A\ is generic. The topology of the strata 
themselves can be as complicated as arbitrary affine varieties over 
\Q\, even for matroids of rank three, by a celebrated theorem of 
Mn\"ev \cite{Mnev}. These strata are connected by ``deletion maps," whose fibers 
are themselves complements of arrangements \cite{BB,F2}.

Realizations in $\Gamma(G)$ correspond to arrangements which have the 
same underlying matroid $G$, as determined by the arbitrary ordering 
of 
the hyperplanes.
Thus, for the study of homotopy type as a function of intrinsic 
combinatorial structure (i.e., without regard to labelling), the true 
``moduli space" for arrangements should be the quotient of 
$\G_\ell(\C^n)$ by the action of the $S_n\times(\C^*)^n$. Then linear 
isomorphism classes of arrangements with isomorphic underlying 
matroids (or isomorphic intersection lattices) correspond to points 
of the orbit space $\Gamma(G)/\Aut(G)$.

Randell's lattice-isotopy theorem can be reformulated as follows: two 
arrangements which are connected by a path in $\widehat{\Gamma}(G)$ 
(or $\Gamma(G)$) 
have diffeomorphic complements. Thus one is led to the difficult 
problem of understanding the set of path components of 
$\Gamma(G)/\Aut(G)$.

More detailed combinatorial data will suffice to uniquely determine 
the homotopy type of the complement. For instance, in the case of 
complexified real arrangements, the defining forms $\alpha_i, 1 \leq 
i \leq n$ determine an underlying {\em oriented matroid}. This is 
most easily described in terms of bases: the matroid $G(\A)$ is 
determined by the collection \B\ of maximal independent subsets $B 
\subseteq [n]$. These can naturally be identified with ordered 
subsets of $[n]$. The {\em oriented} matroid $\widehat{G}(\A)$ is 
then a 
partition $\B=\B_+\cup \B_-$ of the set of ordered bases of $G(\A)$ 
into positive and negative bases, corresponding to the sign of the 
(nonzero) determinant of the corresponding ordered sets of linear 
forms. The work of Salvetti \cite{S1}, as refined by Gelfand and 
Rybnikov \cite{GR}, shows that the underlying oriented matroid of a 
complexified real arrangement uniquely determines the homotopy type 
of the complement. In fact one can construct a partially ordered set 
${\mathcal K}(\widehat{G})$ directly from the oriented matroid 
$\widehat{G}$ 
whose ``nerve", or collection of linearly ordered subsets, forms a 
simplical complex homotopy equivalent to the complement. In 
subsequent work, Bj\"orner and Ziegler \cite{BZ} (see also Orlik 
\cite{O2}) 
generalized the 
construction to arbitrary arrangements (or arrangements of 
subspaces), in 
terms of combinatorial structures called {\em 2-matroids} \cite{BZ} or {\em 
complex oriented matroids} \cite{Z1}. They showed that this detailed 
combinatorial data determines the complement up to piecewise-linear 
homeomorphism. 

The relation between Randell's lattice-isotopy theorem and the 
combinatorial complexes of \cite{S1,GR,BZ,O2} has not been fully 
explored. In particular, it would be interesting to cast the notion 
of lattice-isotopy in combinatorial terms, i.e., as a sequence of 
elementary ``isotopy moves" on the posets $\mathcal{K}(\widehat{G})$ 
which leave the homotopy type of the nerve unchanged. A first step in 
this direction was accomplished in \cite{F3}. We pose this as our 
first open problem.

\begin{problem} Prove a combinatorial lattice-isotopy theorem, that 
``isotopic" (complex) oriented matroids (with the same underlying matroid) 
determine homotopy equivalent cell complexes.
\end{problem}

\end{subsection} 

\begin{subsection}{Homotopy classification}
\label{htpy}
The fundamental question whether the homotopy type of $M(\A)$ is 
uniquely determined by $G(\A)$ was answered in the negative by 
Rybnikov in \cite{Ryb}. The basic building block of his 
construction is the MacLane matroid, whose affine diagram is pictured 
in Figure \ref{mac}. 
\begin{center}
\begin{figure}[h]
\epsfig{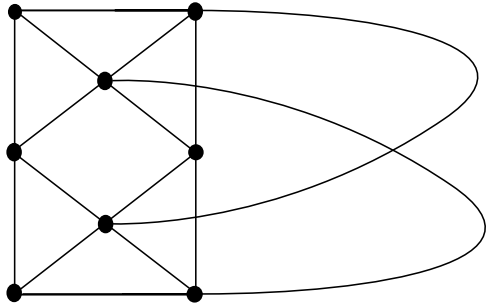}
\caption{The MacLane matroid}
\label{mac}
\end{figure}
\end{center}
For this matroid $G$, the realization space $\R(G)$ consists 
of two conjugate complex realizations $R$ 
and $\overline{R}$, corresponding to arrangements \A\ and 
$\overline{\A}$. One can ``amalgamate'' these realizations along one 
of the 
three-point lines (rank-two flats) to form arrangements $\A \ast 
\A$ and $\A \ast 
\overline{\A}$ of rank four with thirteen hyperplanes. These 
arrangements have the same underlying matroid, of rank four on 13 
points, pictured in Figure \ref{ryb}. 
\begin{center}
\begin{figure}[ht]
\epsfig{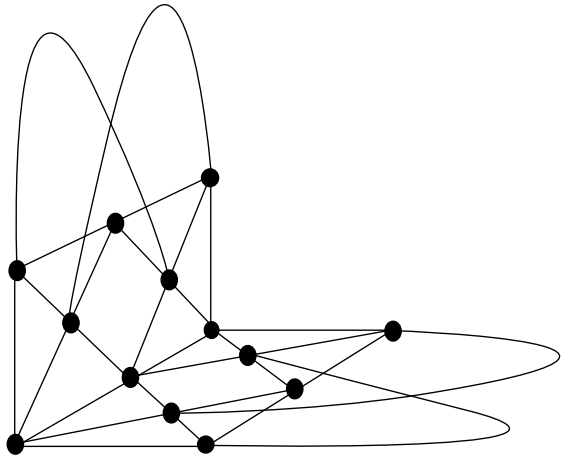}
\caption{The Rybnikov matroid}
\label{ryb}
\end{figure}
\end{center}
Rybnikov establishes some special properties of this matroid, for 
instance, 
that any 
automorphism of the $OS$ algebra arises from a matroid automorphism, 
which must preserve or interchange the factors of the amalgamation. 
Using these he is able to show that the arrangements $\A \ast \A$ and 
$\A \ast 
\overline{\A}$ have nonisomorphic fundamental groups, since the first 
has an automorphism which switches the factors preserving 
orientations of the natural generators, while the only automorphism 
of the 
second which switches factors must reverse orientations. Refer to 
Section \ref{pie1} for a more detailed description of the fundamental 
group. Rybnikov actually uses the rank-three truncation of this 
matroid, and 3-dimensional generic sections of these arrangements, 
but this operation does not affect the fundamental group.

The last part of Rybnikov's argument is quite delicate and very 
specialized. None of the known invariants of 
fundamental groups, for instance those described elsewhere in this 
paper, will distinguish these two groups.

\begin{problem} Find a general invariant of arrangement groups that 
distinguishes the two Rybnikov arrangements, and generalize his 
construction.
\end{problem}

To date this is the only known example of this phenomenon. In 
particular it is not known if this behavior is exhibited by 
complexified arrangements. 

\begin{problem} Prove that the underlying matroid of a complexified 
arrangement determines the homotopy type, or find a counter-example.
\end{problem}

Partial results along these lines were obtained by Jiang and Yau \cite{JY2}
and Cordovil \cite{C4}. In \cite{JY2} a condition on the underlying matroid 
$G$ is given which implies that the realization space of $G$ is path-connected, 
so that any two arrangements realizing $G$ have diffeomorphic 
complements by the 
lattice-isotopy theorem. In \cite{C4} it is shown that complexified 
arrangements whose underlying matroids are isomorphic via a 
correspondence which 
preserves a (geometrically defined) ``shelling order" will have 
identical braid-monodromy groups.

The extent to which arrangements with non-isomorphic matroids can 
have 
homotopy equivalent complements has also been studied 
(see, e.g., \cite{F6,F3,CS2,EF,F7}) with some degree of success. One 
approach to this problem is purely combinatorial, namely to classify 
$OS$ algebras up to graded algebra isomorphism. This approach is 
adopted in \cite{F6,F7,EF}. A powerful invariant is developed in 
\cite{F7}, sufficient to distinguish all known non-trivial examples which 
are not 
already known to be isomorphic.

At this point all known examples of matroids with 
isomorphic $OS$ algebras can be explained by two simple operations 
\cite{FP,Pender}. The first 
of these is a construction involving a well-known equivalence of 
affine arrangements 
arising from the ``cone-decone'' construction \cite[Prop. 5.1]{OT}, 
along 
with the trivial fact that the complement of the direct sum of affine 
arrangements, denoted $\coprod$ in \cite{OT}, is diffeomorphic to the 
cartesian product of the complements of the factors. In fact this 
construction can be applied to arbitrary pairs of matroids to yield 
central arrangements with non-isomorphic matroids and diffeomorphic 
complements \cite{EF,FP}. This construction always yields 
arrangements with non-connected (i.e., nontrivial direct sum) 
matroids. Jiang and Yau \cite{JY1} show that this phenomenon cannot occur in 
rank three, that is, the diffeomorphism type of the complement of a rank-three 
arrangement uniquely determines the underlying matroid. Thus the rank-three 
examples of \cite{F3}, which have non-isomorphic underlying matroids, have 
complements which are homotopy equivalent but not diffeomorphic.

The second operation which yields isomorphic $OS$ algebras is 
truncation. It is shown in \cite{Pender} that the truncations of two 
matroids 
with isomorphic $OS$ algebras will have the same property. (It is not known if 
truncation preserves homotopy equivalence). These two 
``moves'' suffice to explain the examples produced in \cite{OT,F3}, 
indeed all known examples of this phenomenon. Thus it seems an 
orderly 
classification of $OS$ algebras may be within reach.
\begin{problem} Classify $OS$ algebras up to graded isomorphism.
\end{problem}

\noindent
In the alternative, we suggest the following.
\begin{problem} Find a pair of arrangements with homotopy equivalent 
complements  and whose underlying matroids are non-isomorphic, 
connected, and inerectible (i.e., not truncations).
\end{problem}

Cohen and Suciu in \cite{CS,CS2,CS3} approach this same problem of 
homotopy classification using invariants of the fundamental group. 
Their approach has the advantage that it may also be used to 
distinguish the complements of arrangements with the same underlying 
matroid. Some of this work is described elsewhere in this paper.
Here we merely remark on the surprising connection described in 
\cite{CS3,L5,L6} 
between the 
characteristic 
varieties of \cite{L4} arising from the Alexander invariant of the fundamental group, and the 
resonant varieties of \cite{F7}, which arise from the $OS$ algebra.

\end{subsection}

\end{section}

\begin{section}{Algebraic properties of the group of an arrangement}

The topology of hyperplane complements seems to be to a large extent 
controlled by the fundamental group. These ``arrangement groups" have 
relatively simple global structure, being pieced together out of free 
groups in a fairly straightforward way (see Sections \ref{pie1} and 
\ref{weight}), but 
have surprisingly delicate fine structure. At the time of the writing 
of \cite{FR1} there was a great deal of activity around the 
study of the lower central series of these groups, and connections 
with rational homotopy theory and Chen's theory of iterated 
integrals. In this section we report on progress in these areas in 
the intervening years. 

\begin{subsection}{The LCS formula, quadratic algebras, rational 
$K(\pi ,1)$ and parallel arrangements}
\label{lcs}
Discoveries of Kohno \cite{Ko} and the authors \cite{FR2} showed that 
Witt's formula for the lower central series of finitely generated 
free Lie algebras (or, equivalently, free groups) generalized to a 
wide class of hyperplane complements. 
The so-called LCS formula reads $$\prod_{n\geq 1} 
(1-t^n)^{\phi_n}=\sum_{i\geq 0} b_i(-t)^i,$$
relating the ranks $\phi_n$ of factors in the lower central series of 
the fundamental group $\pi_1(M)$ to the betti numbers 
$b_i=\dim(A^i(\A))$ of $M$. In \cite{FR2,J2} it is shown that this 
formula holds for all fiber-type arrangements. These are arrangements 
whose 
underlying matroids are supersolvable \cite{T2}. This result was 
ostensibly 
extended to {\em rational $K(\pi,1)$} arrangements in \cite{F4,Ko2}. 
(See also Section \ref{dnrefl}.) We refer the reader to \cite{F4,OT} for a 
precise definition of 
rational $K(\pi,1)$ arrangement. Briefly, if \S\ is the 1-minimal 
model of $M$ (or, equivalently, of $A(\A)$), then \A\ is rational 
$K(\pi,1)$ if $H^*(\S)\cong 
A(\A)$. It is shown in \cite{F4} that fiber-type arrangements 
are rational $K(\pi,1)$.

The technical results of \cite{FR2} were used in \cite{FR5} to show 
that fundamental groups of fiber-type arrangements (in particular, 
the pure braid group) are residually nilpotent. This result turned 
out to be important for the theory of knot invariants of finite type 
\cite{Stanf}.

The situation surrounding the LCS formula was very much in flux 
during the preparation of \cite{FR1}, a fact reflected in the 
equivocal 
footnotes in the table of implications in that paper. The situation 
has been clarified somewhat in the meantime. Our purpose here is to 
briefly summarize the current understanding of these issues.

Recall that an arrangement of rank three is {\em parallel} if 
for any four hyperplanes of \A\ in general position, there is a fifth 
hyperplane in \A\ containing two of the six pairwise intersections.
The $OS$ algebra $A(\A)$ is {\em quadratic} if the 
relation ideal $I$ (defined in Section \ref{top}) is generated by its elements 
of degree two.
We will sometimes say \A\ is quadratic. This is a combinatorial 
condition, which will be discussed in further 
detail in Section~\ref{quadratic}. In general the quotient of the 
exterior algebra
$\Lambda(e_1,\ldots, e_n)$ by the ideal generated by the degree two 
elements of $I$ is called the {\em quadratic closure} of $A(\A)$, 
denoted $\overline{A}(\A)$. Here is a summary of cogent results 
established in \cite{F4,F8}.

\begin{enumerate}
\label{falk}
\item If \A\ is a rational $K(\pi,1)$ arrangement, then \A\ is 
quadratic. 
\item Every parallel arrangement is quadratic.
\item Every rational $K(\pi,1)$ arrangement satisfies the LCS formula.
\item Every quadratic arrangement satisfies the LCS formula at least 
to third degree.
\end{enumerate} 

In \cite{FR1} we cited an unpublished note which claimed that every parallel
arrangement is a rational $K(\pi,1)$.
Using the construction of 
\cite{F4}, in 1994 Falk wrote a {\em Mathematica} program 
to compute $\phi_4$, and checked the smallest example of a parallel, 
non-fiber-type arrangement of rank 3. This arrangement, labelled 
$X_2$ in \cite{FR1}, consists of the planes $x\pm z=0, y\pm z=0, 
x+y\pm 2z=0$, and $z=0$, and is pictured in Figure \ref{kohno}.
We obtained the result $\phi_4=15$, whereas 
the LCS formula would predict $\phi_4=10$. 
\begin{center}
\begin{figure}[h]
\epsfig{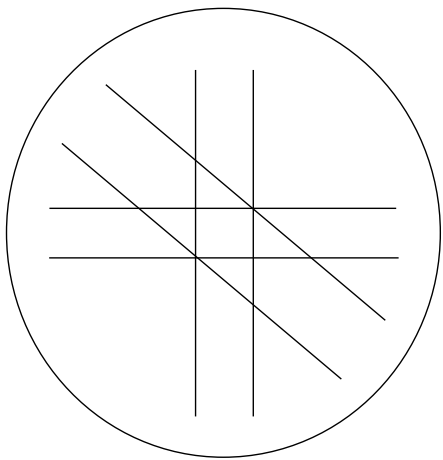}
\caption{The arrangement $X_2$}
\label{kohno}
\end{figure}
\end{center}

So the implications $$\text{parallel} \implies \text{rational 
$K(\pi,1)$},$$
$$\text{quadratic} \implies \text{rational $K(\pi,1)$},$$
$$\text{parallel} \implies \text{LCS},$$
and $$\text{quadratic} \implies \text{LCS}$$
recorded in \cite{FR1}
are all false.

Subsequently, work of Shelton-Yuzvinsky \cite{Y4}, and 
Papadima-Yuzvinsky \cite{Y5} provided further clarification. Let 
${\mathcal L}$ 
denote the holonomy Lie algebra of $M$, the quotient of the free Lie 
algebra on generators $x_1, \ldots x_n$ by the image of the map 
$H_1(M) \to \Lambda^2(H_1(M))$ dual to the cup product. Let $U=U(\A)$ 
be its universal enveloping algebra, a dual object to the 1-minimal 
model \S.  The Hilbert series of $U$ is $\prod_{n\geq 
1}(1-t^n)^{-\phi_n}$.
Kohno constructs a chain complex 
$(R,\delta)$ which, when exact, forms a resolution of $\Q$ as a trivial 
$U$-module.  In this case \A\ is a rational $K(\pi,1)$ 
arrangement, and the LCS formula holds.

Shelton and Yuzvinsky \cite{Y4} realized that $U(\A)$ is the 
Koszul dual of the quadratic closure of $A(\A)$. We refer the reader 
to \cite{Y4} for a precise definition; loosely speaking, the defining 
relations for the Koszul dual $U$ form the orthogonal complement to 
those of $\overline{A}(\A)$ inside the tensor product $T_2(A^1(\A))$.
They observed that the Aomoto-Kohno complex $(R,\delta)$ is the usual 
Koszul complex of $U$, and thus is exact if and only if $U$ is a 
Koszul 
algebra \ --- \ $U$ is {\em Koszul} iff $\Ext_U^{p.q}(\Q,\Q)=0$ unless 
$p=q$.  It follows from this that $A(\A)$ is a quadratic algebra. 
(This observation was also made by Hain \cite{H2}.) 
The LCS formula is then a consequence of Koszul duality.  They give a 
combinatorial proof that $A(\A)$ is quadratic and that $U(\A)$ is Koszul 
if \A\ 
is a supersolvable arrangement.
 
The results of \cite{Y4} were strengthened and 
extended in \cite{Y5} to give a description of $H^*(\S)$ in terms of 
Koszul algebra theory, for more general spaces.  In particular, it 
is shown in \cite{Y5} that \A\ is rational $K(\pi,1)$ if only if 
the 
$OS$ algebra is Koszul.
In addition, Papadima and Yuzvinsky gave an alternate 
proof that the arrangement $X_2$ above fails the LCS formula.  
Finally, using a ``central-to-affine" reduction argument, they were 
able to prove the following.

\begin{theorem} \cite{Y5} For arrangements of rank three, the LCS 
formula holds if and only if the arrangement is fiber-type.
\end{theorem}
Peeva \cite{Pe} applies techniques of commutative algebra and Gr\"obner basis theory to obtain a short proof that supersolvable arrangements satisfy the LCS formula, in addition to other related computational results.

In research closely related to the lower central series of 
arrangement groups, Kohno used the iterated integral/holonomy Lie 
algebra approach to construct representations of the (pure) braid 
group, and more generally to study the monodromy of local systems over 
hyperplane complements.  This work is also closely tied to the theory 
of 
generalized hypergeometric functions.  See \cite{Ko3} for a 
description of these developments.  Cohen and Suciu pursued similar 
ideas using methods more closely connected to those of \cite{FR2} in 
\cite{CS5}.

\end{subsection} 

\begin{subsection}{The $D_n$ reflection arrangements}
\label{dnrefl}
The fundamental 
groups of the reflection arrangements of type $D_n$ have been studied 
using some of the technical machinery of \cite{FR2}.  Note that these 
arrangements, for $n>3$, are not supersolvable.  The author of \cite{Mark} 
constructs a presentation which he claims presents these 
fundamental groups as ``almost direct products" in the sense of 
\cite{FR2,CS5}.  He used this to show that these groups are 
residually 
nilpotent.  In 1994 we tried to use this presentation to get more 
precise calculations for the lower central series of these groups, at 
least for $n=4$.  In fact we found that the presentation in 
\cite{Mark} is not correct.  Even for the $D_3$ arrangement, which is 
supersolvable, the results one deduces from \cite{Mark} do not jibe 
with the LCS formula, which is known to hold for $D_3$. In \cite{Mark2} Liebman 
and 
Markushevich adopt a different approach and derive a different 
presentation to show that the $D_n$ arrangement groups are residually 
nilpotent.

It was in the course of 
this research that we started computing $\phi_4$ by machine.  In 
addition to finding the counterexample $X_2$ described above, we also 
computed $\phi_4=183$ for the $D_4$ reflection arrangement.  The LCS 
formula yields $\phi_4=186$.  So the $D_4$ arrangement fails the LCS 
formula, contrary to another assertion \cite{Ko4} reported on in 
\cite{FR1}. 

The work of Shelton and Yuzvinsky \cite{Y4} make it clear why the argument of 
\cite{Ko4} 
for the LCS formula for the $D_n$ reflection arrangements 
fails: these arrangements, for $n>3$, do not have quadratic $OS$ 
algebras, by \cite{F4}.  Hence the Aomoto-Kohno complex $R_\cdot$ 
cannot be exact for these arrangements.

So we are left with no examples of arrangements which are not 
supersolvable, yet are rational $K(\pi,1)$, and no examples of 
arrangements satisfying the LCS formula which are not rational 
$K(\pi,1)$.

\begin{problem} Find examples of non-supersolvable or non-rational 
$K(\pi,1)$ arrangements 
satisfying the LCS formula, or prove that such examples do not exist.
\end{problem}

\end{subsection} 

\begin{subsection}{Work of Cohen and Suciu on the Chen groups}

As noted above, the ranks of the quotients in the lower central 
series of
fiber-type arrangements are determined by the betti numbers of the
complement. From this point of view, the pure braid groups look like
products of free groups (though they are not; see \cite{FR5}.) In the 
last
few years, Cohen and Suciu have introduced the Chen groups into the 
study of
arrangements, providing a computable tool for distinguishing similar
arrangements.

The Chen groups of a group $G$ are the lower central series quotients 
of $G$
modulo its second commutator subgroup $G^{\prime \prime }$. If for 
any group 
$G$ we let $\Gamma _{k}(G)$ denote the $k^{th}$ lower central series
subgroup, then the homomorphism $G\rightarrow G/G^{\prime \prime }$ 
induces
an epimorphism 
\begin{equation*}
\frac{\Gamma _{k}(G)}{\Gamma _{k+1}(G)}\rightarrow \frac{\Gamma
_{k}(G/G^{\prime \prime })}{\Gamma _{k+1}(G/G^{\prime \prime 
})}=k^{th}\text{
Chen group}
\end{equation*}

Thus the ranks $\phi_k$ of quotients of lower central series groups are no 
less than
the corresponding ranks $\theta _{k}$ of Chen groups. In the case of 
the
pure braid group, the ranks $\theta _{k}$ are determined in 
\cite{CS}; they
are given by the generating function 
\begin{equation*}
\sum_{k=2}^{\infty }\theta _{k}t^{k-2}=\binom{n+1}{4}\cdot 
\frac{1}{(1-t)^{2}%
}-\binom{n}{4}
\end{equation*}

In particular, these numbers differ from those for the product of free
groups, providing a tidy proof that the pure braid groups are not such
products.

Cohen and Suciu \cite{CS4} provide a detailed study of these groups 
including
a method for their computation from a presentation of the Alexander
invariant (see the discussion of presentations of the fundamental 
group
below.) It is interesting that while these groups are very effective 
in
distinguishing similar groups, there is not yet an example of
combinatorially equivalent arrangements with different Chen ranks. In
particular, they do not distinguish the examples of Rybnikov \cite{Ryb} of
combinatorially equivalent, homotopically different arrangements (see Section 
\ref{htpy}).

\end{subsection} 

\begin{subsection}{Cohomological properties of the fundamental group}
\label{cohom}

In 1972 Deligne \cite{D} proved that for a complexification of a real
simplicial arrangement, the complement $M$ is aspherical (also 
expressed by
saying that $M$ is a $K(\pi ,1)$ space.) That 
is,
the universal cover of $M$ is contractible. Since all real reflection
arrangements are simplicial, this solved a question raised and 
partially
answered by Brieskorn in \cite{B}. The original study of this sort of
problem was the work of Fadell and Neuwirth \cite{FN} on the pure 
braid
group. Following \cite{T1}, the authors introduced in \cite{FR2} the notion of 
fiber-type
arrangement and observed that for this class $M$ is aspherical, 
essentially
by the iterated fibration argument of Fadell and Neuwirth. So it is 
natural
to ask: for what arrangements is $M$ aspherical? It is known by work 
of
Hattori \cite{Ha} that not all are \ --- \ the arrangement defined by 
$Q=xyz(x+y+z)$ is the simplest example.

Here we wish to touch upon the algebraic consequences of asphericity. 
Now if 
$M$ is aspherical, the (known) cohomology of $M$ is isomorphic to the
cohomology of the group. Since $M$ has cohomological dimension $\rk(\A) 
<\infty, \pi _{1}(M)$ does also. In addition, $\pi _{1}(M)$ has no torsion, 
and there is a $K(\pi,1)$ space, $\pi=\pi_1(M)$, 
with the 
homotopy type of a finite complex (namely, $M$).
So here is another open
problem:

\begin{problem}
Are all arrangement groups torsion-free?
\end{problem}

The answer is of course yes for real reflection arrangements and for
fiber-type (or supersolvable) arrangements. One approach to this 
question is to show 
that all
arrangement groups are orderable. Here we say a group $G$ is \emph{orderable}
provided that there is a linear order $<$ on $G$ so that \(g<h\) implies
\(cg<ch\) for all \(c \in G\).  It follows easily that an orderable
group has no torsion. 
The braid group was shown orderable by Dehornoy in \cite{Deh}; at the Tokyo 
meeting L. Paris
proved that the group of a fiber-type arrangement is orderable \cite{P4}.  It is 
not known whether all
arrangement groups are orderable.
Note
that the group of an arrangement has a finite presentation of a fairly 
restricted type, 
as
described in Section \ref{pie1}, and that the relators all lie in the
commutator subgroup.

There are some useful observations concerning these ideas in 
\cite{R2}. For instance, we have the following theorem.

\begin{theorem} For $j\geq 2$ the Hurewicz map
\begin{equation*}
\phi :\pi _{j}(M)\rightarrow H_{j}(M)
\end{equation*}
is trivial. 
\end{theorem}

As a consequence, the second homology of  $\pi _{1}(M)$ is
isomorphic to $H_{2}(M)$. In addition, it is mentioned there that the
arrangement defined by $$Q=xyz(y+z)(x-z)(2x+y)$$ has the property 
that
there is no arrangement with aspherical complement with the same intersection 
lattice 
in
rank one and two. The following result is also proved in \cite{R2}.

\begin{theorem} The complement of a central arrangement of rank three is 
aspherical provided that the fundamental group has cohomological dimension three 
and is of type FL.
\label{rand}
\end{theorem}

A group $\pi$ is {\em type FL} provided that $\Z$ (as a trivial
$\Z[\pi]$-module) has a finite resolution by free $\Z[\pi]$-modules. An 
equivalent
statement is that there should exist a {\em finite} CW complex which is a
$K(\pi,1)$-space.  Theorem \ref{rand} shows that for central rank three 
arrangements
asphericity is determined by the fundamental group.
\end{subsection}
\end{section}

\begin{section}{Arrangements with aspherical complements}

Much of the early history of the topology of arrangements revolves 
around the ``$K(\pi,1)$ problem," the problem of determining which 
arrangements have aspherical complements. (Such an arrangement is 
called a {\em $K(\pi,1)$ arrangement}.) This history is described 
in some detail in \cite{FR1} (see also Section \ref{cohom}).  
In addition, we proved an {\em ad hoc} 
necessary condition \cite[Thm.~3.1]{FR1} for asphericity 
involving ``simple triangles," 
and introduced the notion of formal arrangement, which was shown to 
be 
a necessary condition for $K(\pi,1)$ and rational $K(\pi,1)$ 
arrangements.
A great deal of progress was made in these areas in the intervening 
years,
which we report on in this section.

\begin{subsection}{Free arrangements are not aspherical}

In our earlier survey, we highlighted the Saito conjecture, that all 
free
arrangements are aspherical. In 1995 Edelman and Reiner \cite{ER} 
provided
counterexamples, which we briefly describe.

Let $S$ denote the polynomial ring of $V.$ A linear map $\theta
:S\rightarrow S$ is a derivation if for $f,g\in S,$ we have $\theta
(fg)=f\theta (g)+g\theta (f)$. The module of $\A-$derivations is
defined by 
\begin{equation*}
D(\A)=\{\theta \mid \theta (Q)\in QS\}
\end{equation*}
where $Q$ is the defining polynomial of the arrangement. Then the
arrangement is \emph{free }provided that $D(\A)$ is a free 
$S$-module.

It is known \cite{T1} that reflection arrangements are free; for 
their many
pleasant properties see \cite{OT}. In 1975 K. Saito conjectured that 
free
arrangements should be aspherical. In their study of tilings of 
centrally
symmetric octagons in \cite{ER}, Edelman and Reiner found the 
family of arrangements given by $$Q(\A_\alpha)=xyz(x-y)(x-z)(y-z)(x-\alpha 
y)(x-\alpha
z)(y-\alpha z)$$ with $\alpha \in \Re.$ They proved that the corresponding 
arrangements are free for 
all $\alpha $, while they are not aspherical for $\alpha \neq -1,0,1$. The 
proof
of freeness is direct, using addition-deletion \cite[Theorem 
4.51]{OT} while
the non-asphericity follows from the ``simple triangle'' criterion of 
\cite{FR1}. The counter-example $\A_{-2}$ is pictured in Figure \ref{edel}.

\begin{center}
\begin{figure}[h]
\epsfig{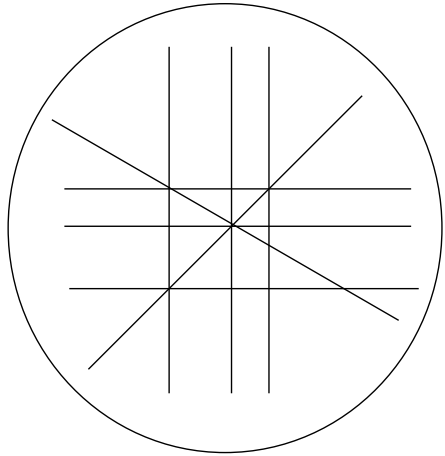}
\caption{Free but not $K(\pi,1)$}
\label{edel}
\end{figure}
\end{center}

\end{subsection}

\begin{subsection}{Formality and related concepts}
\label{quadratic}
The fundamental group of arrangement is determined by a generic 
3-dimensional section.  Based on the idea that $K(\pi,1)$ 
arrangements 
should be extremal in some sense, we developed the notion of formal 
arrangement in \cite{FR1}.  This has been the subject of several 
papers since \cite{BB,BT,Y1,F12}, which provide a 
better understanding of the concept.  Here is a ``modern" 
definition, 
equivalent to the original from \cite{FR1}.

Let $\Phi: \C^n \to V^*$ be given by $\Phi(x)=\sum_{i=1}^n x_i\alpha_i$, 
where the $\alpha_i$ are the defining forms for \A.  Let 
$K=\ker(\Phi)$ and let $F$ be the subspace of $K$ spanned by its 
elements of weight three (i.e., having three nonzero entries).
Then the arrangement \A\ is {\em formal} if $F=K$.

The orthogonal complement $K^\perp \subseteq \C^n$ coincides with the 
point $P_\A \in \G_\ell(\C^n)$ defined in Section \ref{top}. 
Thus the arrangement \A\ is isomorphic to the arrangement in 
$K^\perp$ formed by the coordinate hyperplanes. In the same way, the 
orthogonal complement $F^\perp \supseteq K^\perp$ defines an 
arrangement $\A_F$, called the {\em formalization} of \A. So \A\ is 
formal if and only if $\A=\A_F$. If \A\ is not formal, $\A_F$ 
has strictly greater rank, and \A\ is a (not necessarily generic) 
section of $\A_F$. Also, \A\ and $\A_F$ have isomorphic generic 
``planar" (i.e., rank-three) sections. 

These properties of formalization were asserted 
in \cite{FR1}, but the arguments we had in mind were not correct. The 
clarification described here is due to Yuzvinsky \cite{Y1}.
Examples in \cite{Rad} show that non-formal arrangements need 
not be generic sections of their formalizations. The 
arrangement of Example 2.19 of \cite{Rad} has the property that the {\em free} 
erection of the underlying matroid is not realizable, but (contrary to the 
assertion in \cite{Rad}) there 
is nevertheless a realizable (formal) erection. Matroid 
``erection'' is the reverse of (corank one) truncation; truncation is the 
matroid-theoretic analogue of generic section. The free erection of an erectible 
matroid is the unique erection with ``the most general position" \ --- \ see 
\cite{Wh1}.

These observations are enough to establish the following results from 
\cite{FR1}. The third assertion follows immediately from the second.

\begin{enumerate}
\item If \A\ is a $K(\pi,1)$ arrangement, then \A\ is 
formal.
\item If \A\ is quadratic, then \A\ is formal.
\item If \A\ is a rational $K(\pi,1)$ arrangement, then \A\ is 
formal.
\end{enumerate} 

We asked whether free arrangements are also necessarily formal. This 
was 
established by Yuzvinsky.

\begin{theorem} \cite{Y1} If \A\ is a free arrangement, then \A\ is 
formal.
\label{Yuz}
\end{theorem}

The preceding result was generalized by Brandt and Terao \cite{BT}. 
They define 
the notion of $k$-formal arrangement.
A formal arrangement has the property that all relations among the 
defining equations are consequences of relations which are 
``localized" at rank-two flats, in the sense that an element of $K$ of 
weight three gives rise to a three-element subset of a rank-two 
flat. A formal arrangement is {\em $3$-formal} if all relations among 
these 
local generators of $F=K$ are themselves consequences of relations 
which are localized at rank-three flats of \A.  This construction is 
iterated to define the notion of $k$-formal arrangement for every 
$k\geq 2$. See \cite{BT} for the precise definition. An arrangement 
of rank $r$ is automatically $k$-formal for every $k\geq r$.
The original notion of formality coincides with the case $k=2$.

\begin{theorem} \cite{BT} If \A\ is a free arrangement of rank $r$, 
then \A\ is 
$k$-formal for every $2\leq k < r$.
\end{theorem}
\noindent The converse is false \cite{BT}.

Related work appears in \cite{BB}, where the authors 
show that the discriminantal arrangements  of Manin and Schechtman 
\cite{MS} (see Section \ref{discrim}) are formal, and the ``very generic" 
discriminantal arrangements are $3$-formal, though none are free.

An arrangement is {\em locally formal} \cite{Y1} if, for every 
flat $X \subseteq 
[n],$ the arrangement $\A_X=\{H_i \ | \ i \in X\}$ is formal. Since 
freeness, quadraticity, and $K(\pi,1)$-ness are all ``hereditary 
properties," in that they are inherited by the localizations $\A_X$, 
one has that every free, quadratic, or $K(\pi,1)$ arrangement is 
locally formal.

We asked in \cite{FR1} whether formality is a ``combinatorial 
property'', depending only on the underlying matroid. Yuzvinsky 
constructed counter-examples in \cite{Y1}.

\begin{theorem} \cite{Y1} There exist arrangements $\A_1$ and $\A_2$ with the 
same underlying matroid, such that  $\A_1$ is formal and $\A_2$ is 
not formal.
\end{theorem}

In Figure \ref{yuzfig} are the dual point configurations of Yuzvinsky's 
arrangements. The dotted line in Figure \ref{yuzfig}(b) indicates where to 
``fold" the configuration to erect it to a rank-four configuration. The 
nontrivial planes in the erection are $$12389, \ 12456, \ 13458, \ 13678, \ 
14579, \ 23567, \ 24789, \ 25689, \ \text{and} \ 34679.$$ Note that these two 
configurations are lattice-isotopic (over \C), so neither is free or $K(\pi,1)$.

\begin{center}
\begin{figure}[h]
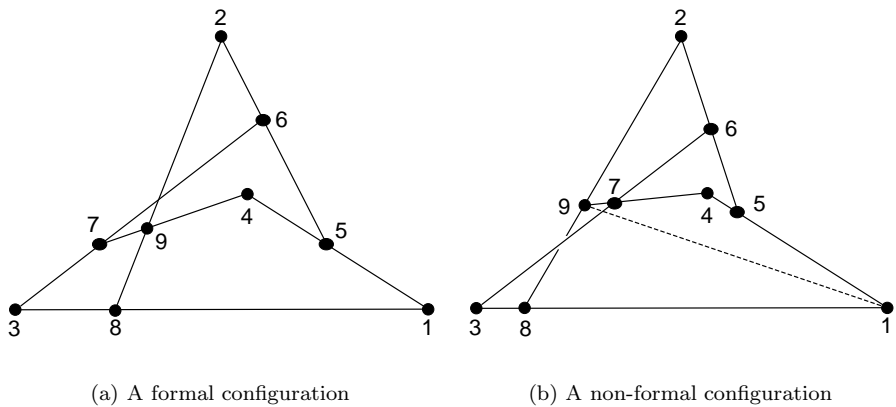

\mbox{\subfigure[A formal configuration]{\epsfig{file=yuz1.epsf, 
height=4.5cm}}\quad
\subfigure[A non-formal configuration]{\epsfig{file=yuz2.epsf, height=4.5cm}}}
\caption{Formality is not matroidal}
\label{yuzfig}
\end{figure}
\end{center}

If \A\ is not formal, then the underlying matroid of \A\ is a {\em 
strong map image} (under the identity map) of that of $\A_F$ (see 
\cite{Ox} for the general definition), 
and the two matroids have the same 
rank-three truncations. These combinatorial properties gave rise to 
several attempts 
to replace the notion of formality with some clearly matroidal 
condition, and strengthen Theorem \ref{Yuz} and assertion (i) above. 
For example one can ask for conditions on a matroid $G$ so that every 
(complex) realization of $G$ is formal. One is naturally led to the notion of 
line-closure.

Let $G$ be a matroid on ground set $[n]$. The {\em 
line-closure} of a subset $S$ of $[n]$ is the smallest subset of 
$[n]$ 
which contains every line (that is, rank-two flat) spanned by points of $S$. A 
set is 
line-closed if it is equal to its line-closure. The matroid $G$ is 
{\em line-closed} if every line-closed subset of $[n]$ is a flat of 
$G$. 
In his current work in progress \cite{F12}, the first author has established the 
following result.

\begin{theorem} An arrangement \A\ is quadratic 
only if the underlying matroid $G(\A)$ is line-closed. 
\label{lineclosed}
\end{theorem}

\begin{corollary} The underlying matroid of a rational $K(\pi,1)$ 
arrangement is necessarily line-closed.
\end{corollary}

The converse of Theorem \ref{lineclosed}, that \A\ is quadratic when $G(\A)$ is 
line-closed, is very likely also true. A crucial step in the proof is yet to be 
completed, however, so this assertion remains an open problem.

Yuzvinsky \cite{Y6} defined a {\em formal matroid} to be a matroid 
$G$ 
possessing a basis
(of $\rk(G)$ points) whose line-closure is $[n]$. Every line-closed 
matroid is formal in this sense. In fact a matroid $G$ is line-closed 
if and 
only 
if the line-closure of {\em every} basis of each flat $X$ is equal to $X$. 
Every realization of a formal matroid is formal.

In \cite{F12} we define a matroid $G$ to be {\em taut} if $G$ is not 
a strong 
map 
image of a matroid $G'$ of greater rank with the same points and 
lines, and {\em locally taut} if every flat of $G$ is taut. 
Every line-closed matroid is locally taut, in fact every formal 
matroid is taut. Every realization of a 
(locally) taut 
matroid is (locally) formal. There 
exist matroids which are taut but not formal \cite{Cra}. A weak version of the 
first part of the 
following problem was suggested by Yuzvinsky in his talk \cite{Y6}. 

\begin{problem} Prove that the matroid of a free or $K(\pi,1)$ 
arrangement is necessarily taut.
\end{problem}

Joseph Kung has pointed out to us that a locally taut matroid is 
uniquely determined by its points and lines, which suggests the 
following interesting problem.

\begin{problem} Prove that the underlying matroid of a locally formal 
arrangement (e.g. a free or $K(\pi,1)$ arrangement) is uniquely 
determined by its points and lines.
\end{problem} 

This last problem is a variant on the following questions from \cite{FR1}, 
the first of which is Terao's Conjecture, and both of which remain 
open.

\begin{problem} Prove that freeness and $K(\pi,1)$-ness of 
arrangements are matroidal properties.
\end{problem}

We will refrain from discussing Terao's Conjecture further, except to 
pose a weak version which fits the spirit of this paper, and is 
interesting in its own right.

\begin{problem} Prove that freeness is preserved under lattice-isotopy.
\end{problem}
\end{subsection}

\begin{subsection}{Tests for asphericity}
\label{weight}
Some progress was also made on the problem of finding sufficient 
conditions for an arrangement to be $K(\pi,1)$. The main results are 
the weight test of \cite{F1} and its application to factored 
arrangements by Paris \cite{P2}. A new technique involving modular 
flats was recently discovered and presented at the conference 
\cite{P3,FP}.

The complement $M$ of 
a 2-dimensional affine arrangement \A\ is built up out of $K(\pi,1)$ 
spaces, specifically $(r,r)$ torus link complements, in a relatively 
simple way, 
as is reflected in the Randell-Salvetti-Arvola presentations (see 
Section \ref{pie1}). In fact this structure mirrors precisely  
constructions from geometric group theory related to complexes of 
groups. This observation allows one to construct a relatively 
well-behaved cell complex which has the homotopy type of the 
universal 
cover of $M$, and to apply the weight test of Gersten and Stallings 
\cite{St} to derive a test for asphericity of $M$. 

\begin{theorem} \cite{F1} If \A\ is a complexified affine 
arrangement in $\C^2$ that admits an \A -admissible, aspherical 
system of 
weights, then \A\ is a $K(\pi,1)$ arrangement.
\end{theorem}

The question remains what an \A -admissible, aspherical system of 
weights is. This involves the complex $B$ of bounded faces in 
the subdivision of $\Re^2$ determined by \A. A weight system is an 
assignment of a real 
number weight to each ``corner'' of each 2-cell in $B$. The 
system is aspherical if the sum of weights around any $d$-gon at most $d-2$. 
The system is \A -admissible if certain sums of weights 
at 
vertices of $\Gamma$ are at least $2\pi$. See \cite{F1} for more 
detail.

The universal cover complex constructed in \cite{F1} may be used in 
some cases to 
construct explicit essential spheres showing that $M$ is not 
aspherical. Radloff \cite{Rad} used this method to prove some 
necessary conditions for $K(\pi,1)$-ness, along the lines of 
the ``simple triangle'' test of \cite{FR1}, and found several new examples of 
non-$K(\pi,1)$ 
arrangements.

Falk and Jambu introduced the notion of {\em factored} arrangement in 
\cite{FJ}, originally in an attempt to find a combinatorial criterion 
for freeness. A {\em factorization} of an arrangement \A\ is a 
partition of $[n]$ such that each flat of $G(\A)$ of rank $p$ meets 
precisely $p$ 
blocks, and meets one of them in a singleton, for each $p$. This 
property is necessary and sufficient for the $OS$ algebra $A(\A)$ to 
have a complete tensor product factorization - see 
\cite{BZ2,FJ,T4,OT}.
When \A\ has a factorization, we say \A\ is {\em factored}.

Paris realized that a factorization of a rank-three arrangement 
provides a 
template for a very simple \A -admissible, aspherical weight system.

\begin{theorem} \cite{P2} If \A\ is a factored, complexified arrangement in 
$\C^3$, then \A\ is a $K(\pi,1)$ arrangement.
\end{theorem}

Every supersolvable arrangement is factored, so this result provides 
a new, wider class of $K(\pi,1)$ arrangements, at least in rank three.

\begin{problem} Show that factored arrangements of arbitrary rank are 
$K(\pi,1)$.
\end{problem}

A flat $X$ of a matroid $G$ is {\em modular} if $\rk(X\vee 
Y)+\rk(X\wedge Y)=\rk(X)+\rk(Y)$ for every flat $Y$. The following result was 
discovered independently by Paris and Falk-Proudfoot

\begin{theorem} \cite{P3,Proud,FP} If $X$ is a modular flat of arbitrary 
rank in $G(\A)$, then there is a topological fibration $M(\A) \to M(\A_X)$ whose 
fiber is the complement of a projective arrangement.
\end{theorem}

This generalizes the corank-one case, which gives rise to fiber-type 
arrangements, 
established in \cite{T2}. The new result can be used to construct or recognize 
$K(\pi,1)$ arrangements if the base (whose matroid is the modular 
flat $X$)
and fiber (whose matroid is the complete principal truncation of $G(\A)$
along $X$) are known to be $K(\pi,1)$. This method is used to 
construct some interesting new examples in \cite{FP}. Refer to Paris' 
paper \cite{P4} in this volume for more details.

\end{subsection}

\begin{subsection}{Some crucial examples}

In this section we want to briefly discuss some specific and 
interesting
types of arrangements for which the $K(\pi,1)$ problem is unsolved. 
These might be regarded as test subjects for
new techniques; they qualify as ``the first unknown cases.'' 

First we 
cite another improvement to the table of implications in \cite{FR1}. 
Recall the definition of parallel arrangement from Section \ref{lcs}. 
In \cite{FR1} we had listed the implication ``parallel $\implies$ $K(\pi,1)$'' 
as 
``not known, of significant interest.'' In unpublished work, 
Luis Paris has shown this implication
to be false. Specifically, he showed that the Kohno arrangement $X_2$ 
(defined in Section \ref{lcs}) is not $K(\pi,1)$. The proof establishes 
that the fundamental group contains a subgroup isomorphic to ${\mathbb 
Z}^4$; the result then follows from \cite[Thm.~3.2]{FR1}. The copy of ${\mathbb 
Z}^4$ is generated by $a,b,c,$ and the commutator $[d,e]$,
where $a,b,c,d,$ and $e$ are the canonical generators corresponding to 
the hyperplanes $x\pm z=0, \ z=0,$ and $x+y\pm 2z=0$ respectively.

\subsubsection{Complex reflection arrangements}

Fadell and Neuwirth showed in 1962 that the 
complement
of the  $A_{\ell }$ reflection (or braid) arrangement is $K(\pi ,1)$. In 1973 
Brieskorn
proved this for many real reflection arrangements, followed soon 
thereafter
by Deligne's proof of the general case. Orlik and Solomon extensively
studied arrangements of hyperplanes invariant under finite groups generated by 
complex
reflections (see \cite[Chapter 6]{OT}). It is natural to ask if all 
such
arrangements are aspherical. We believe the conjecture that they are 
is due to Orlik, though it was proposed long before it ever appeared 
in print. It is known \cite{OT} that the answer is
affirmative in all cases except six exceptional, non-complexified arrangements, 
some of which have rank three. The
proofs for the known cases use a variety of techniques, and 
essentially
proceed from the Shepard-Todd classification of irreducible unitary reflection 
groups (see, e.g., \cite{OT}). What seems to be 
missing is a
unifying property, similar to the simplicial property for real 
reflection
arrangements exploited by Deligne. The closest approach to this goal is the 
work
reported in \cite[p. 265]{OT} which proves the asphericity of 
arrangements
associated to Shephard groups (symmetry group of a regular convex 
polytope.)
Here the problem is reduced to the (already solved) problem for an
associated real reflection arrangement.

\begin{problem} Give a uniform proof that all unitary reflection 
arrangements are $K(\pi,1)$.
\end{problem}

\end{subsection}

\subsubsection{Discriminantal arrangements}

\label{discrim}
Experience seems to show us that questions involving asphericity are 
quite
complex for all arrangements but tractible for restricted classes
(reflection, fiber-type, generic). One interesting class is that of 
the
discriminantal arrangements introduced by Manin and Schechtmann 
\cite{MS}.
Rather than give the full definition here we will describe the rank three
examples, where the problem is already interesting. 

Consider a real affine arrangement of lines in the plane, obtained by
taking a collection of $n$ points, no three of which are collinear, 
and
drawing all $\binom{n}{2}$ lines through pairs of these points. Then 
embed
this configuration in the plane $z=1$ in three-space and cone over the
origin to obtain a central real three-arrangement. Then 
complexify.

This process can result in arrangements with distinct matroidal and topological 
structure, even for fixed $n$ \cite{F2,BB}. The discriminantal arrangements are 
obtained from ``very generic" collections of points, for which no three of the 
$\binom{n}{2}$ lines are concurrent except at the original $n$ points.

The arrangement $C(4)$ is linearly equivalent to the braid arrangement of rank 
three. An easy calculation shows that the 
Poincar\'{e} polynomial associated to the cohomology of $C(n)$ does not 
factor over $\mathbb Z$ for $n\geq 5$, so that
these arrangements are not free and are not of fiber-type. Also $C(n)$ is not 
simplicial for $n\geq 5$.
The arrangements $C(n)$ for $n\geq 6$ are not aspherical, by 
\cite[Thm.~3.1]{FR1}.

For $n=5$, one obtains a complexified central 
three-arrangement of $10$ planes. This arrangement is not factored. More 
generally $C(5)$ does not support an admissible, aspherical system of weights, 
so the weight test fails. On the other hand, all of the standard necessary 
conditions for asphericity hold. 

\begin{problem} Determine whether the discriminantal arrangement $C(5)$ is 
$K(\pi,1)$.
\end{problem}

A  solution to this problem would also determine whether the  space of 
configurations of six points in general position in $\C P^2$ is aspherical 
\cite{F2}, a result which would be of significant interest.

\subsubsection{Deformations of reflection arrangements}

A ``deformation" of a reflection arrangement is an affine arrangement with 
defining equations of the form $$\alpha_i(x_1,\ldots, x_\ell)=c_{ij},$$ where 
the $\alpha_i$ are the positive roots in some root system, and $c_{ij}\in \Re$. 
This class of arrangements is of great interest to combinatorialists, and is the 
subject of the paper of Athanasiadis in this volume \cite{Athan}. 

As is our custom, we ``cone" to obtain a central arrangement. For instance, 
based on the root system of type $B_2$, we obtain the {\em $B_2$ Shi 
arrangement}, defined by
the  polynomial $$Q=xyz(x+y)(x-y)(x-z)(y-z)(x+y-z)(x-y-z).$$  (Shi arrangements 
are obtained by setting $c_{i1}=0$ and $c_{i2}=1$ for all $i$.) This nine-line 
complexified
arrangement has a factorization, given by the partition 
$$\{\{4\},\{1,2,5,7\},\{3,6,8,9\}\},$$ and is therefore a $K(\pi,1)$ arrangement. 
On the other hand, the Shi arrangement constructed in a similar way from the 
root system of type $G_2$ is not factored or simplicial, and has no simple 
triangle.

\begin{problem} Decide whether the $G_2$ Shi arrangement is $K(\pi,1)$.
\end{problem}

\noindent More generally, we propose the following.

\begin{problem} Decide which Shi arrangements are $K(\pi,1)$.
\end{problem}
\end{section} 

\begin{section}{Topological properties of the group of an arrangement}

At the time of the publication of \cite{FR1}, a presentation of the fundamental 
group of the complement of a complexified arrangement had been derived 
\cite{R3}. In the meantime, a similar presentation was found for arbitrary 
complex arrangements \cite{Ar3}, and several different ``spines" for the 
complement, some of them modelled on group presentations, were constructed 
\cite{S1,F3,CS2,L1}. These group presentations have been used to study the Milnor 
fibration and Alexander invariants of the complement. We report briefly on these 
ideas here.

\begin{subsection}{Presentations of $\pi _{1}$}
\label{pie1}
We have seen earlier in the discussion of the lower central series, 
Chen
groups and group cohomology that certain classes of arrangements
(fiber-type, simplicial) have well-behaved fundamental groups. Due to 
work
of Arvola \cite{Ar3}, Randell \cite{R3} and Salvetti \cite{S1} an 
explicit
presentation of $\pi_{1}(M)$ can be written. See \cite[Section 5.3]{OT} for a 
clear
exposition of Arvola's presentation for any complex arrangement, and 
\cite{F3} for the explicit presentation and some applications of Randell's
presentation, which holds for complexified arrangements and is naturally 
simpler
than the general case. A different approach, using the notion of ``labyrinth," 
is adopted by Dung and Vui in \cite{VD} to arrive at  similar presentations for 
arbitrary arrangement groups.

In these presentations one first takes a planar section (or, more precisely, the 
projective image), so that one 
is
working with an affine arrangement in $\C^{2}$. Then there is one 
generator
for each line of the arrangement, and one set of relations 
for
each intersection. In all cases the relations consist entirely of
commutators, but to date this has not shed much light on the 
questions of
group cohomology, torsion in the fundamental group, or other 
properties
(such as orderability) of the fundamental group. A general theme for
questions is: to what extent do arrangement groups mimic the 
properties of
the pure braid groups. 

The concept of braid monodromy was introduced by B.~Moishezon \cite{Moish}. Libgober showed in \cite{L1} that the braid monodromy presentation of the fundamental group yields a two-complex with the homotopy type of the complement of an algebraic curve (e.g., a line arrangement) transverse to the line at infinity.

Motivated in part by \cite{L1}, the first author showed in \cite{F3} that for arbitrary line arrangements the 2-complex modelled on the presentation of 
\cite{R3} serves 
as an
efficient model for constructing the homotopy type of the complement 
(in the
case of 3-arrangements). This construction was then used to 
construct a
number of examples with different intersection lattice but same 
homotopy
type (see also Section \ref{htpy}).

In related work Cohen and Suciu \cite{CS2} have given an explicit
description of the braid monodromy of a complex arrangement, using Hansen's 
theory of polynomial covering maps. They show 
that the resulting presentation of the fundamental group is 
equivalent to the Randell-Arvola presentation via Tietze transformations that do 
not affect the homotopy type of the associated 2-complex. It follows that the 
complement is
homotopy equivalent to the 2-complex modelled on either of these presentations, 
generalizing the result of \cite{F3}.
For this work Cohen and Suciu used 
extensively the concept of {\em braided wiring diagram}, which
we briefly describe below.  The notion of braided wiring diagram generalizes 
Goodman's concept 
of wiring diagram \cite{Good},
and was earlier considered for arrangements in
\cite{C3}. (Wiring diagrams appear in combinatorics as geometric models for 
rank-three oriented matroids.)  The presentations
of \cite{R3} and \cite{Ar3} use versions of this idea.  In brief, the braided 
wiring
diagram can be thought of as a template for the fundamental group (or, for line
arrangements, the homotopy type.)

Here is a sketch of the construction. For examples and further details, in 
particular, a  beautiful derivation using polynomial covering space theory, see
\cite{CS2}. Since we are interested in the 
fundamental group,
consider an affine arrangement \A\ in $\C^{2}$. Choose coordinates in $\C^{2}$ 
so that the 
projection to the first coordinate is generic.  Suppose that the images 
$y_{1}, \ldots,y_{n}$ of the 
intersections of the lines have distinct real parts.  Choose a basepoint 
$y_{0}\in \C \setminus \{y_{1}, \ldots,y_{n} \}$, and assume the real parts of 
$y_{i}$
are decreasing with $i$.  Let $\xi$ be a smooth path which begins with $y_{0}$ 
and passes
in order through the $y_{i}$, horizontal near each $y_{i}$. Then the braided 
wiring 
diagram is $\mathcal{W}=\{(x,z) \in \xi \times \C \mid Q(x,z)=0\}.$ (Recall that 
$Q$
is the defining polynomial of the arrangement.) 

This braided wiring diagram should be viewed as a picture of the braid monodromy
of the fundamental group of the arrangement (or as a picture of the fundamental 
group itself).
In a sense, it carries the attaching (or amalgamating)
information as one computes the fundamental group using the Seifert-Van Kampen 
theorem.  Each
actual node in the wiring diagram gives a set of relators, as does each 
crossing.
In particular, it is shown in \cite{CS2} that the braided wiring diagram 
recovers the Arvola or Randell 
presentation of $\pi_{1}(M)$.  Indeed, in the real case, the braided wiring 
diagram 
can be identified with the usual drawing of the arrangement in $\Re^{2}$.

As is the case with ordinary braids,
there are ``Markov moves'' with which one can modify such a wiring diagram 
to realize any braid-equivalence of the underlying braid monodromies.  
These are given explicitly in \cite{CS2}. Rudimentary moves of this type, called ``flips," first appeared in \cite{F3}. Among 
the consequences we
note the following results which relate braid monodromy and braided wiring 
diagrams
to lattice isotopy of line arrangements (that is, arrangements in $\C^{2}$).

\begin{theorem} \cite{CS2} Lattice-isotopic arrangements in $\C^{2}$ 
have braided wiring diagrams which
are related by a finite sequence of Markov moves and their inverses.
\end{theorem}

\begin{theorem} \cite{CS2} Line arrangements with braid-equivalent monodromies
have isomorphic underlying matroids.
\end{theorem}

\end{subsection} 

\begin{subsection}{The Milnor fiber}

The defining polynomial $Q=\prod_{i=1}^n \alpha _{i}$ is homogeneous of degree 
$n$ and can be considered 
as a
map 
\begin{equation*}
Q:M\rightarrow \C^{*}
\end{equation*}
It is well-known that this map is the projection of a fiber bundle, 
called
the Milnor fibration, and that the Milnor fiber $F=Q^{-1}(1)$ should 
be of
interest. In \cite{R5} it was shown that this Milnor fibration is 
constant
in a lattice-isotopic family, so that the Milnor fiber is indeed an
invariant of lattice-isotopy. Because of this we propose the
following definition, analogous to the definition made in the theory 
of
knots.

\begin{definition}
Two arrangements are called \emph{(topologically) equivalent} if they 
are
lattice-isotopic. We say the arrangements have the same 
\emph{(topological)
type.}
\end{definition}

Thus, arrangements are topologically equivalent if and only if 
they lie in the same path 
component of
some matroid stratum in the Grassmannian. With this terminology, we have the 
following result.

\begin{theorem} \cite{R5} The Milnor fiber and fibration are invariants of
topological type.
\end{theorem}

Now, $F$ is simply an $n-$fold cover of the complement of the 
projectivized
arrangement in $\C P^{\ell -1}.$ Since the algorithms of the previous 
section
work to compute the fundamental group of this latter space, questions
involving the fundamental group and cohomology of $F$ are also 
questions
involving the group of the arrangement. In particular, while the 
cohomology
of $M$ is determined by the intersection lattice, that of $F$ may not 
be.
The situation is analogous to that of plane curves, where work going 
back to
Zariski \cite{Zar1} shows that not only the type but the position of 
the
singularities affects the irregularity. (The irregularity here is 
simply
half the ``excess'' in the first betti number of $F$.)

Early results concerning the Milnor fiber of an arrangement (often in the general context of plane curves) appear in work of Libgober \cite{L1,L2,L3,L4} and Randell \cite{Ran}, particularly with respect to Alexander invariants. Libgober's work gave considerable information about the homology of the Milnor fiber in relation to the number and type of singularities of the arrangement, their position and the number of lines. The paper \cite{Ran} observed that the Alexander polynomial was equal to the characteristic polynomial of the monodromy on the Milnor fiber.

The paper of Artal-Bartolo \cite{AB1} included an interesting example: for the rank three braid arrangement $A_3$ the first betti number of the Milnor fiber is seven, an excess of two over the five ``predicted" by the number of lines. (This result can be obtained as an interesting exercise by applying the Reidemeister-Schreier rewriting algorithm to the presentations of the fundamental group.) Orlik and Randell \cite{OR} showed that in the generic case the cohomology of the Milnor fiber is minimal, given by the number of lines, below the middle dimension.

Cohen and Suciu carry forward the study of the Milnor fiber in \cite{CS4}. 
Using the
group presentation and methods of Fox calculus they give 
twisted
chain complexes whose homology gives that of the Milnor fiber. Their 
methods
are effective, and several explicit examples are given. The monodromy 
action
on the Milnor fiber is of course crucial, 
and
this monodromy is determined as well.

Finally, we note the following problem, which remains open after many 
years.

\begin{problem} Prove that the homology of the Milnor fiber of \A\ 
depends only on the underlying matroid.
\end{problem}

\end{subsection}

\end{section}

\begin{ack}
\em
The idea to hold a birthday conference in honor of Peter Orlik was initially 
suggested by the second author in 1995. We would like to thank Mutsuo Oka And 
Hiroaki Terao for their hard work in organizing the meeting. We also wish to 
thank the referee for his helpful observations concerning deformations of 
arrangements. 

The two authors were both students of Peter Orlik in Madison, Randell in the 
early 1970's and Falk in the early 1980's. We are happy to have the chance to 
thank him, in print, for introducing us to the field of hyperplane arrangements 
and for his enthusiasm, friendship and support through the years.
\end{ack}

\end{document}